\documentclass[a4paper,10pt]{article}
\usepackage[english]{babel}
\usepackage{mlmodern}
\usepackage[T1]{fontenc}
\usepackage{amsmath}
\usepackage{amsfonts}
\usepackage{titlesec}
\usepackage{amssymb,fancyhdr}
\usepackage{amsthm,booktabs}
\usepackage{multirow,mathtools,xfrac}
\usepackage{placeins}
\usepackage{graphicx,caption}
\usepackage{subcaption}
\usepackage{epstopdf}
\usepackage{setspace}
\usepackage[nottoc]{tocbibind}

\numberwithin{equation}{section}
\numberwithin{figure}{section}
\numberwithin{table}{section}

\usepackage{abstract}

\usepackage{paralist}
\usepackage{babel,csquotes}
\usepackage{lipsum}
\usepackage{url}
\usepackage{etoolbox}
\patchcmd{\abstract}{\null\vfil}{}{}{}
\usepackage[format=hang,font=small,labelfont=bf]{caption}
\usepackage[left=1in, right=1in, top=1in, bottom=0.5in, includefoot]{geometry}
\usepackage[colorlinks=true,allcolors=blue,citecolor=blue,linkcolor=blue]{hyperref}
\usepackage{authblk} 
\title{\textbf{\Large Galerkin-Bernstein Approximations of the System of Time Dependent Nonlinear Parabolic PDEs }}
       \author{\small Hazrat Ali}
        \author{\small Nilormy Gupta Trisha}
        \author{\small Md. Shafiqul Islam*}
         \affil{Department of Applied Mathematics, University of Dhaka, Dhaka 1000, Bangladesh.}
        \affil{ Corresponding author and Email : mdshafiqul@du.ac.bd}
        \date{}

\makeatletter
\patchcmd{\@maketitle}
  {\addvspace{0.5\baselineskip}\egroup}
  {\addvspace{-1\baselineskip}\egroup}

\makeatother
\begin{document}
\maketitle
		\lhead{\thepage}
\section*{Abstract}
The purpose of the research is to find the numerical solutions to the system of time dependent nonlinear parabolic partial differential equations (PDEs) utilizing the Modified Galerkin Weighted Residual Method (MGWRM) with the help of modified Bernstein polynomials. An approximate solution of the system has been assumed in accordance with the modified Bernstein polynomials. Thereafter, the modified Galerkin method has been applied to the system of nonlinear parabolic PDEs and has transformed the model into a time dependent ordinary differential equations system. Then the system has been converted into the recurrence equations by employing backward difference approximation. However, the iterative calculation is performed by using the Picard Iterative method. A few renowned problems are then solved to test the applicability and efficiency of our proposed scheme. The numerical solutions at different time levels are then displayed numerically in tabular form and graphically by figures. The comparative study is presented along with $\displaystyle L_2$ norm, and $\displaystyle L_\infty$ norm.\\

\noindent{\bf\normalsize Keywords:} {\normalsize Parabolic PDE System, Modified Galerkin Method, Modified Bernstein Polynomial, Backward Difference Method, Gray-Scott Model}
\section{Intrduction}
Reaction-diffusion systems have been extensively studied during the $\displaystyle 20^{th}$ century. The study of the reaction-diffusion system reveals that different species have interactions with one another and that after these interactions, new species are created via chemical reactions. The solution of the reaction-diffusion system shows the chemical reaction's underlying mechanism and the various spatial patterns of the chemicals involved. 
Animal coats and skin coloration have been linked to reaction-diffusion processes, which have been considered to constitute a fundamental basis for processes associated with morphogenesis in biology. 
\\
\noindent There are numerous notable examples of coupled reaction-diffusion systems such as the Brusselator model, Glycolysis model, Schnackenberg model, Gray-Scott model, etc. With the help of the system size expansion, a stochastic Brusselator model has been suggested and investigated in the study cited in \cite{r39}. The reaction-diffusion Brusselator model has been addressed by Wazwaz et al. through the decomposition technique \cite{r40}. Because of its potential to provide a close analytical solution, the fractional-order Brusselator model was studied by Faiz et al \cite{r41}. The Brusselator system stability of a reaction-diffusion cell as well as the Hopf bifurcation analysis of the system have been detailed by Alfifi \cite{r42}. Qamar has analyzed the dynamics of the discrete-time Brusselator model with the help of the Euler forward and nonstandard difference schemes \cite{r43}. The research article cited in \cite{r36} has been prepared by investigating the numerical analysis of the Glycolysis model using a well-known finite difference scheme.  
Adel et al \cite{r45} have examined the synchronization problem of the Glycolysis reaction-diffusion model and designed a novel convenient control law. David et al \cite{r47} have analyzed the stability of turing patterns of the Schnackenberg model. Liu et al \cite{r48} have developed the bifurcation analysis of the aforementioned model. Khan et al. \cite{r49} have established a scheme for the solution of the fractional order Schnackenberg reaction-diffusion system. Numerical explorations have been applied to analyze the pattern formations of the model in the research article cited in \cite{r50}. Gray and Scott \cite{r7} were the first to introduce the Gray-Scott model. They have proposed this model as an alternative to the autocatalytic model of Glycolysis \cite{r8}. For this model, Pearson \cite{r9} has employed experimental studies to depict several sophisticated spot-type structures. Mazin et al. \cite{r10} have conducted an experiment using a computer simulation to investigate a range of far-from-equilibrium occurrences that emerge in a bistable Gray-Scott model. Many renowned authors \cite{r11,r15} have evaluated the preceding model in which self-replicating structures have been noticed. McGough et al. \cite{r13} have conducted research on the bifurcation analysis of the patterns that are depicted in the model. In the research cited in \cite{r21}, the linear stability and periodic stationary solutions of this model have been investigated. Some analytical results of this model have also been explored \cite{r17}. Several prominent authors have studied the spatiotemporal chaos of the model in the research studies cited in \cite{r18} and \cite{r19}. Furthermore, Wei \cite{r20} has analyzed the pattern formation of the two-dimensional Gray-Scott model. The model has also been explored by Kai et al. \cite{r22} using an innovative technique known as the second-order explicit implicit methodology. In recent years, the nonlinear Galerkin finite element approach has become increasingly prevalent as a means to investigate the model \cite{r2,r1}. Mach \cite{r5} has performed an in-depth examination of the quantitative evaluation of the model's numerical solution. In references \cite{r3} and \cite{r12}, the Gray-Scott reaction-diffusion system has been the subject of extensive wave modeling studies by eminent scholars. The simulation of the coupled model has been carried out by Owolabi et al. \cite{r14} using the higher-order Runge-Kutta method. The well-known Gray-Scott model's numerical findings have been calculated using the help of the hyperbolic B-spline \cite{r4}. In order to analyze the ionic version of the model while it is being affected by an electric field, the Galerkin method has been deployed \cite{r6}. With the use of the hybrid-asymptotic numerical method, Chen et al. \cite{r24} have investigated the model's dynamic behavior and stability. In the research study cited in \cite{r23}, special polynomials have been employed to numerically solve the Gray-Scott model. Han et al. \cite{r25} have conducted an exhaustive investigation on the three-dimensional Gray-Scott model. In the process of assessing the model, the cubic B-spline has proven to be of considerable use by Mittal et al \cite{r26}.\\

\noindent In the disciplines of engineering and mathematical physics, the Weighted Residual Method is an approximation method that can be leveraged to resolve problems. Analysis of structures, thermal expansion, stream of fluids, movement of masses, and the electromagnetic potential, etc. are examples of prominent problem fields of concern. Several distinct Weighted Residual Method variations are within our reach. The Galerkin Weighted Residual Method (also known as GWRM) has been put into practice for centuries, long before the invention of computers. It is generally agreed that this strategy is one of the best and most often used approaches available. Lewis and Ward have provided a comprehensive overview of the process in the article that is referenced in \cite{r27}. This methodology has been effectively implemented in the well-known Black-Scholes model by Hossan et al. \cite{r28}. Shirin et al. \cite{r29} have employed the Galerkin method in conjunction with other special polynomials to analyze the Fredholm equations. In the research referred to in \cite{r30}, the approach was utilized to solve boundary value problems. In addition, this method has been used to perform a numerical calculation of the eigenvalues associated with the Sturm-Liouville problem \cite{r31}. There have been several successful uses of this method for problems involving metal beams and polygonal ducts with rounded edges \cite{r32, r33}.\\
The objective of this study is to employ the modified Galerkin Weighted Residual Method in conjunction with the appropriate special polynomials to numerically evaluate the one-dimensional reaction-diffusion systems. Based on our best information, this study is presently unavailable. In addition to that, the study has provided the validation necessary to use the approach in one-dimensional reaction-diffusion systems. The main merit and advantage of the study are that by solving this type of system of equations, we will be able to analyze the behavior of the ecological system and forecast its future. 
\\
\noindent The article is split up into four sections. Section 2 provides a detailed explanation of the formulation of our proposed method to solve the system of nonlinear parabolic partial differential equations. In the third section, the approach's implications are shown while analyzing the aforementioned system. Numerical and graphical representations are included here as well. The fourth section contains some concluding remarks and a general discussion.
\section{Mathematical Formulation}
\noindent Let us commence with the following system over the domain $[-L, L]$
\begin{align}
\setstretch{2.0}
\left.
\begin{array}{lcl}
\displaystyle \frac{\partial M}{\partial t}=\varepsilon_1\frac{\partial^2 M}{\partial x^2} - f(M,N) + p(1-M)\\
\displaystyle \frac{\partial N}{\partial t}=\varepsilon_2\frac{\partial^2 N}{\partial x^2}+ f(M,N) -(p+q)N
\end{array}
\right\}
\label{eq:3.1}
\end{align}
\noindent The boundary and initial conditions are as follows:
\begin{align}
\left.
\begin{array}{lcl}
\setstretch{1.0}
    \displaystyle M(-L,t)=M(L,t)=\theta_0\\
    \displaystyle N(-L,t)=N(L,t)=\gamma_0
    \end{array}
    \right\}
    \label{eq:3.2}
    \end{align}
    and
    \begin{align}
    \setstretch{1.0}
    \left.
    \begin{array}{lcl}
    \displaystyle M(x,0)=M_0(x)\\
    \displaystyle N(x,0)=N_0(x)
    \end{array}
    \right\}
    \label{eq:3.3}
\end{align}
Let us assume the approximate solutions of System (\ref{eq:3.1}) be of the form
\begin{align}
\left.
\begin{array}{lcl}
\setstretch{0.9}
    \displaystyle \widetilde{M}(x,t)=\theta_0+\sum_{j=0}^{n}c_j(t)B_j(x)\\
    \displaystyle \widetilde{N}(x,t)=\gamma_0+\sum_{j=0}^{n}d_j(t)B_j(x)\\
    \end{array}
    \right\}
    \label{eq:3.4}
    \end{align}
    \noindent where $\displaystyle B_j$'s are the modified Bernstein polynomials and $\displaystyle c_j$ and $\displaystyle d_j$ are the coefficients dependent on time. 
    The first terms of the approximate solutions (\ref{eq:3.4}) have come from the boundary conditions of the system. The modified Bernstein polynomials are defined as follows:
    \begin{align}
    B_{n,m}(x)= \binom{m}{n} \frac{(x-L)^n(U-x)^{m-n}(x-L)(U-x)}{(U-L)^m}\hspace{1cm} n=0,1,2,..., m\nonumber
   \end{align}
   \FloatBarrier
where $U$ \& $L $ are the upper and lower limits of $x$.
    The last terms of Solution (\ref{eq:3.4}) will vanish at the boundary points.
    Therefore, the residual functions are
    \begin{align}
\setstretch{1.8}
\left.
\begin{array}{lcl}
\displaystyle R_1(x,t)=\frac{\partial \widetilde{M}}{\partial t}-\varepsilon_1\frac{\partial^2 \widetilde{M}}{\partial x^2} + f(\widetilde{M},\widetilde{N}) - p(1-\widetilde{M})\\
\displaystyle R_2(x,t)=\frac{\partial \widetilde{N}}{\partial t}-\varepsilon_2\frac{\partial^2 \widetilde{N}}{\partial x^2}- f(\widetilde{M},\widetilde{N}) +(p+q)\widetilde{N}\\
\end{array}
\right\}
\label{eq:3.5}
\end{align}
\noindent Now we form the residual equations as:
\begin{align}
    &\int_{-L}^{L}R_1(x,t) B_i(x)dx=0\\
    &\int_{-L}^{L}R_2(x,t) B_i(x)dx=0 
    \label{12}
\end{align}
From the first residual equation, we can write
\begin{align}
 &\int_{-L}^{L}\Bigg[\frac{\partial \widetilde{M}}{\partial t}-\varepsilon_1\frac{\partial^2 \widetilde{M}}{\partial x^2} + f(\widetilde{M}, \widetilde{N}) - p(1-\widetilde{M})\Bigg]B_i(x)dx=0
 \label{eq:31.5}
 \end{align}
 \noindent Now we apply integration by parts in the above equation
 \begin{align}
 \int_{-L}^{L}\frac{\partial \widetilde{M}}{\partial t}B_idx+\int_{_L}^{L}\varepsilon_1\frac{\partial \widetilde{M}}{\partial x}\frac{\partial B_i}{\partial x}dx + \int_{-L}^{L} f(\widetilde{M}, \widetilde{N}) B_idx-\int_{-L}^{L}p(1-\widetilde{M})B_idx=\varepsilon_1\Big[\frac{\partial \widetilde{M}}{\partial x}B_i\Big]_{-L}^{L}
    \label{eq:31.6}
\end{align}
Then we substitute solution (\ref{eq:3.4}) in Equation (\ref{eq:31.6}). Therefore, the equation becomes,
\begin{align*}
  & \int_{-L}^{L}\frac{\partial}{\partial t}\Big(\theta_0+\sum_{j=0}^{n}c_jB_j\Big)B_idx+\int_{-L}^{L}\varepsilon_1\frac{\partial}{\partial x}\Big(\theta_0+\sum_{j=0}^{n}c_jB_j\Big)\frac{\partial B_i}{\partial x}dx  +\int_{-L}^{L}f\Big(\theta_0+\sum_{j=0}^{n}c_jB_j, \gamma_0+\sum_{j=0}^{n}d_jB_j\Big)B_i dx\nonumber\\
  &\hspace{5cm}-\int_{-L}^{L}p\Big(1-\big(\theta_0+\sum_{j=0}^{n}c_jB_j\big)\Big)B_idx=\varepsilon_1\Big[\frac{\partial}{\partial x}\Big(\theta_0+\sum_{j=0}^{n}c_jB_j\Big)B_i\Big]_{-L}^{L} \nonumber\\
  & \text{or,} \int_{-L}^{L}\frac{\partial \theta_0}{\partial t}B_idx+\int_{-L}^{L}\sum_{j=0}^{n}\frac{\partial c_j}{\partial t}B_j B_idx+\int_{-L}^{L}\varepsilon_1\frac{\partial \theta_0}{\partial x}\frac{\partial B_i}{\partial x}dx+\sum_{j=0}^{n}c_j\int_{-L}^{L}\varepsilon_1\frac{\partial B_j}{\partial x}\frac{\partial B_i}{\partial x}dx\nonumber\\
     &+\int_{-L}^{L}f\Big(\theta_0+\sum_{j=0}^{n}c_jB_j, \gamma_0+\sum_{j=0}^{n}d_jB_j\Big)B_i dx-\int_{-L}^{L}pB_idx+\int_{-L}^{L}p\theta_0 B_idx+\sum_{j=0}^{n}c_j\int_{-L}^{L}pB_jB_idx\nonumber\\
    &\hspace{9cm}=\varepsilon_1\Bigg[\frac{\partial \theta_0}{\partial x}B_i\Bigg]_{-L}^{L}+\varepsilon_1\Bigg[\sum_{j=0}^{n}c_j\frac{\partial B_j}{\partial x}B_i\Bigg]_{-L}^{L}
\end{align*}
This finally becomes
\begin{align}
&\int_{-L}^{L}\frac{\partial \theta_0}{\partial t}B_idx+\int_{-L}^{L}\sum_{j=0}^{n}\frac{\partial c_j}{\partial t}B_j B_idx+\int_{-L}^{L}\varepsilon_1\frac{\partial \theta_0}{\partial x}\frac{\partial B_i}{\partial x}dx+\sum_{j=0}^{n}c_j\int_{-L}^{L}\varepsilon_1\frac{\partial B_j}{\partial x}\frac{\partial B_i}{\partial x}dx\nonumber\\
     & +\int_{-L}^{L}\Gamma\Big(\theta_0, \gamma_0, \sum_{k=0}^{n}c_kB_k, \sum_{l=0}^{n}d_lB_l\Big)B_i dx+\sum_{j=0}^{n}d_j\int_{-L}^{L}\Omega\Big(\theta_0, \gamma_0, \sum_{k=0}^{n}c_kB_k, \sum_{l=0}^{n}d_lB_l\Big)B_j B_idx\nonumber\\
     &-\int_{-L}^{L}pB_idx+\int_{-L}^{L}p\theta_0 B_idx +\sum_{j=0}^{n}c_j\int_{-L}^{L}pB_jB_idx =\varepsilon_1\Bigg[\frac{\partial \theta_0}{\partial x}B_i\Bigg]_{-L}^{L}+\varepsilon_1\Bigg[\sum_{j=0}^{n}c_j\frac{\partial B_j}{\partial x}B_i\Bigg]_{-L}^{L}   
     \label{eq:31.7}
     \end{align}
     \noindent The first terms on both sides, and third terms on the left-hand side Equation (\ref{eq:31.7}) become zero because of boundary conditions. Therefore, the equation reduces to,
     \begin{align}
    & \sum_{j=0}^{n}\frac{d c_j}{dt}\int_{-L}^{L}B_j B_idx+\sum_{j=0}^{n}c_j\Bigg(\int_{-L}^{L}\varepsilon_1\frac{dB_j}{dx}\frac{d B_i}{dx}dx+\int_{-L}^{L}pB_jB_idx-\varepsilon_1\Big[\frac{d B_j}{dx}B_i\Big]_{-L}^{L}\Bigg)\nonumber\\
    &+\sum_{j=0}^{n}d_j\int_{-L}^{L}\Omega\Big(\theta_0, \gamma_0, \sum_{k=0}^{n}c_kB_k, \sum_{l=0}^{n}d_lB_l\Big)B_j B_idx=-\int_{-L}^{L}\Gamma\Big(\theta_0, \gamma_0, \sum_{k=0}^{n}c_kB_k, \sum_{l=0}^{n}d_lB_l\Big)B_i dx\nonumber\\
&\hspace{8cm} +\int_{-L}^{L}pB_idx -\int_{-L}^{L}p\theta_0 B_idx
    \label{eq:61.1}
\end{align}
The derivative and non-derivative terms of Equation (\ref{eq:61.1}) can be summarized via standard matrix notation as follows:
\begin{align}
    [C_1]\Big\{\frac{dc_j}{dt}\Big\}+[K_1]\{c_j\}+[K_2]\{d_j\}=[F_1]
    \label{eq:3.9}
\end{align}
where
\begin{align}
    &C_{1_{ij}}= \int_{-L}^{L}B_j B_idx\nonumber\\
    &K_{1_{ij}}= \int_{-L}^{L}\varepsilon_1\frac{dB_j}{dx}\frac{d B_i}{dx}dx+\int_{-L}^{L}pB_jB_idx-\varepsilon_1\Big[\frac{d B_j}{dx}B_i\Big]_{-L}^{L}\nonumber\\
    &K_{2_{ij}}= \int_{-L}^{L}\Omega\Big(\theta_0, \gamma_0, \sum_{k=0}^{n}c_kB_k, \sum_{l=0}^{n}d_lB_l\Big)B_j B_idx\nonumber\\
    &F_{1_{i}}= -\int_{-L}^{L}\Gamma\Big(\theta_0, \gamma_0, \sum_{k=0}^{n}c_kB_k, \sum_{l=0}^{n}d_lB_l\Big)B_i dx+\int_{-L}^{L}pB_idx -\int_{-L}^{L}p\theta_0 B_idx\nonumber
\end{align}
\noindent Here, $K_1$ and $K_2$ are $n \times n$ matrices, $C_1$ is $n \times n$ matrix, and $F_1$ is $n \times 1$ matrix. The first two matrices $K_1$ and $K_2$ are called stiffness matrices. The other two matrices $C_1$ and $F_1$ are called forced matrix, and load vector respectively.\\
\noindent Therefore, we apply the backward difference method on the first term of Equation (\ref{eq:3.9}) and rearrange the resulting terms as follows:
\begin{align}
\setstretch{2.0}
     &[C_1]\Big\{\frac{c_j-c_{j-1}}{\Delta t}\Big\}+[K_1]\{c_j\}+[K_2]\{d_j\}=[F_1]\nonumber\\
     & Or, \Big(\frac{1}{\Delta t}[C_1]+[K_1]\Big)\{c_j\}+[K_2]\{d_j\}=\frac{1}{\Delta t}[C_1]\{c_{j-1}\}+[F_1]
     \label{eq:31.10}
\end{align}
\noindent The second residual equation can be written as,
\begin{align*}
&\int_{-L}^{L}\Bigg[\frac{\partial \widetilde{N}}{\partial t}-\varepsilon_2\frac{\partial^2 \widetilde{N}}{\partial x^2}- f(\widetilde{M},\widetilde{N}) +(p+q)\widetilde{N}\Bigg]B_i(x)dx=0
\end{align*}
\noindent After employing integration by parts and then substitution of (\ref{eq:3.4}) reduces the above equation,
\begin{align*}
  & \int_{-L}^{L}\frac{\partial}{\partial t}\Big(\gamma_0+\sum_{j=1}^{n}d
    _jB_j\Big)B_idx+\int_{-L}^{L}\varepsilon_2\frac{\partial}{\partial x}\Big(\gamma_0+\sum_{j=1}^{n}d_jB_j\Big)\frac{\partial B_i}{\partial x}dx+\int_{-L}^{L}(p+q)\Big(\gamma_0+\sum_{j=0}^{n}d_jB_j\Big)B_idx\nonumber\\
    & \hspace{4cm}-\int_{-L}^{L}f\Big(\theta_0+\sum_{j=0}^{n}c_jB_j, \gamma_0+\sum_{j=0}^{n}d_jB_j\Big)B_i dx=\varepsilon_2\Big[\frac{\partial}{\partial x}\Big(\gamma_0+\sum_{j=0}^{n}d_jB_j\Big)B_i\Big]_{-L}^{L}\nonumber
\end{align*}
\begin{align}
    \text{or,}&\int_{-L}^{L}\frac{\partial \gamma_0}{\partial t}B_idx+\int_{-L}^{L}\sum_{j=0}^{n}\frac{\partial d_j}{\partial t}B_j B_idx+\int_{-L}^{L}\varepsilon_2\frac{\partial \gamma_0}{\partial x}\frac{\partial B_i}{\partial x}dx+\sum_{j=0}^{n}d_j\int_{-L}^{L}\varepsilon_2\frac{\partial B_j}{\partial x}\frac{\partial B_i}{\partial x}dx\nonumber\\
     & -\int_{-L}^{L}\Pi\Big(\theta_0, \gamma_0, \sum_{l=0}^{n}d_lB_l\Big)B_i dx-\sum_{j=0}^{n}c_j\int_{-L}^{L}\Phi\Big( \gamma_0, \sum_{l=0}^{n}d_lB_l\Big)B_j B_idx+\sum_{j=0}^{n}d_j\int_{-L}^{L}(p+q)B_jB_idx\nonumber\\
     &\hspace{4cm}=-\int_{-L}^{L}(p+q)\gamma_0B_idx+\varepsilon_2\Bigg[\frac{\partial \gamma_0}{\partial x}B_i\Bigg]_{-L}^{L}+\varepsilon_2\Bigg[\sum_{j=0}^{n}d_j\frac{\partial B_j}{\partial x}B_i\Bigg]_{-L}^{L}
     \label{eq:3.12}
\end{align}
Since the first, and third terms on the left-hand side and the first term on the right-hand side of Equation (\ref{eq:3.12}) become zero, the equation reduces to,
\begin{align}
&\sum_{j=0}^{n}\frac{d d_j}{dt}\int_{-L}^{L}B_j B_idx+\sum_{j=0}^{n}d_j\Bigg(\int_{-L}^{L}\varepsilon_2\frac{dB_j}{dx}\frac{d B_i}{dx}dx+\int_{-L}^{L}(p+q)B_jB_idx-\varepsilon_2\Big[\frac{d B_j}{dx}B_i\Big]_{-L}^{L}\Bigg)\nonumber\\
    &-\sum_{j=0}^{n}c_j\int_{-L}^{L}\Phi\Big(\gamma_0, \sum_{l=0}^{n}d_lB_l\Big)B_j B_idx =\int_{-L}^{L}\Pi\Big(\theta_0, \gamma_0, \sum_{l=0}^{n}d_lB_l\Big)B_i dx -\int_{-L}^{L}(p+q)\gamma_0B_idx 
\label{eq:3.13}
\end{align}
The derivative and non-derivative terms of Equation (\ref{eq:3.13}) can be summarized via standard matrix notation as follows:
\begin{align}
    [C_2]\Big\{\frac{dd_j}{dt}\Big\}+[K_3]\{c_j\}+[K_4]\{d_j\}=[F_2]
    \label{eq:3.14}
\end{align}
where
\begin{align}
    &C_{2_{ij}}= \int_{-L}^{L}B_j B_idx\nonumber\\
     &K_{3_{ij}}= -\int_{-L}^{L}\Phi\Big(\gamma_0, \sum_{l=0}^{n}d_lB_l\Big)B_j B_idx\nonumber\\
    &K_{4_{ij}}= \int_{-L}^{L}\varepsilon_2\frac{dB_j}{dx}\frac{d B_i}{dx}dx+\int_{-L}^{L}(p+q)B_jB_idx-\varepsilon_2\Big[\frac{d B_j}{dx}B_i\Big]_{-L}^{L}\nonumber\\
    &F_{2_{i}}= \int_{-L}^{L}\Pi\Big(\theta_0, \gamma_0, \sum_{l=0}^{n}d_lB_l\Big)B_i dx -\int_{-L}^{L}(p+q)\gamma_0 B_idx\nonumber
\end{align}
\noindent Here, $K_3$ and $K_4$ are $n \times n$ matrices, $C_2$ is $n \times n$ matrix, and $F_2$ is $n \times 1$ matrix. They are called stiffness matrices, forced matrices, and load vectors respectively.\\
\noindent The application of the backward difference method on the first term of Equation (\ref{eq:3.14}) results in the following equation,
\begin{align}
 \Big(\frac{1}{\Delta t}[C_2]+[K_4]\Big)\{d_j\}+[K_3]\{c_j\}=\frac{1}{\Delta t}[C_2]\{d_{j-1}\}+[F_2]
     \label{eq:3.15}
\end{align}
\noindent By assembling Equations (\ref{eq:31.10}) and (\ref{eq:3.15}), we get the following recurrent system,
\begin{align}
\setstretch{2.0}
\left.
    \begin{array}{lcl}
    \displaystyle \Big(\frac{1}{\Delta t}[C_1]+[K_1]\Big)\{c_j\}+[K_2]\{d_j\}=\frac{1}{\Delta t}[C_1]\{c_{j-1}\}+[F_1]\\
    \displaystyle [K_3]\{c_j\}+\Big(\frac{1}{\Delta t}[C_2]+[K_4]\Big)\{d_j\}=\frac{1}{\Delta t}[C_2]\{d_{j-1}\}+[F_2]
    \end{array}
    \right\}
    \label{eq:3.16}
\end{align}
\noindent To calculate the initial values of $\displaystyle c_j$ and $\displaystyle d_j$, the initial conditions are set in Galerkin sense as follows,
\begin{align}
    & \int_{-L}^{L}\widetilde{M}(x,0)B_idx=\int_{-L}^{L}M_0(x)B_idx\nonumber\\
 \text{or,}& \int_{-L}^{L}\Big(\theta_0+\sum_{j=1}^{n}c_j(0)B_j(x)\Big)B_idx=\int_{-L}^{L}M_0(x)B_idx\nonumber\\
\text{equivalently,}& \sum_{j=0}^{n}c_j(0) \int_{-L}^{L}B_jB_idx=\int_{-L}^{L}M_0(x)B_idx-\int_{-L}^{L}\theta_0 B_idx
    \label{eq:3.17}
\end{align}
\noindent and
\begin{align}
    & \int_{-L}^{L}\widetilde{N}(x,0)B_idx=\int_{-L}^{L}N_0(x)B_idx\nonumber\\
 \text{equivalently,} &\int_{-L}^{L}\gamma_0 B_idx+\int_{-L}^{L}\sum_{j=0}^{n}d_j(0) B_jB_idx=\int_{-L}^{L}N_0(x)B_idx\nonumber\\
\text{or,}&\sum_{j=0}^{n}d_j(0) \int_{-L}^{L}B_jB_idx=\int_{-L}^{L}N_0(x)B_idx-\int_{-L}^{L}\gamma_0 B_idx
    \label{eq:3.18}
    \end{align}
    \noindent This process will help us to evaluate the numerical solutions of the nonlinear reaction-diffusion systems.
\section{Numerical Examples and Applications}
\noindent In this section, the previously described approach has been implemented into practice by solving a few examples of practical issues. Our methodology has been shown to be valid after being applied to the first test problem. The aforementioned procedure is then used, with a variety of parameters, to assess the subsequent test problems. The $\displaystyle L_2$ norm and $\displaystyle L_\infty$ norm has been determined by the following expression,
\begin{align*}
&\displaystyle L_2 \ Norm=||M_{\Delta t}-M_{\frac{\Delta t}{2}}||_2\\
&\displaystyle L_\infty \  Norm=||M_{\Delta t}-M_{\frac{\Delta t}{2}}||_\infty
\end{align*}
Where $\Delta t$ is the time increment and $M_{\Delta t}$ is the approximate solution obtained using time increment $\Delta t$.\\

\begin{figure}[ht!]
		\centering
		\begin{minipage}{0.45\textwidth}
			\centering
			\includegraphics[width=\textwidth]{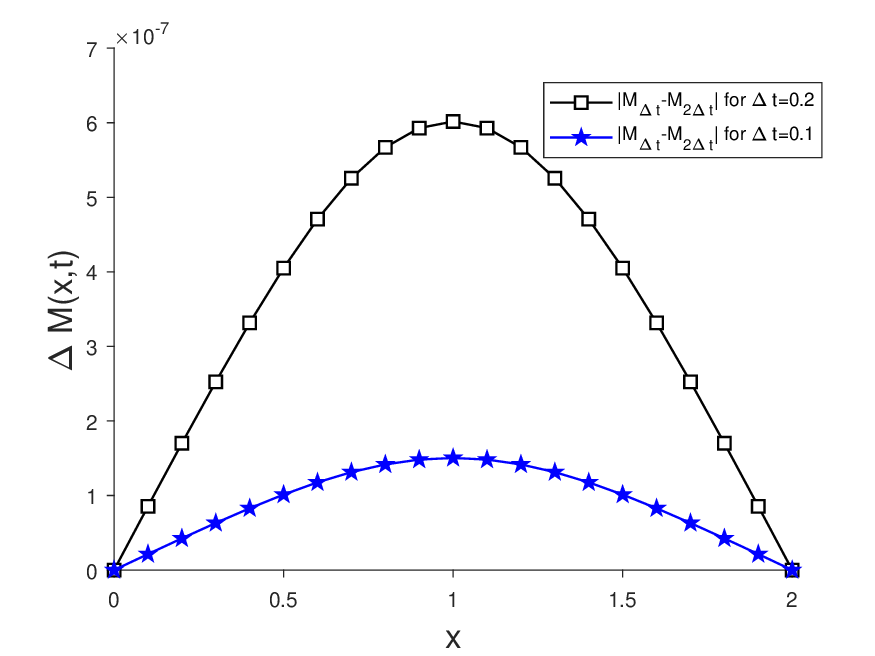}
			\subcaption{}
		\end{minipage}\hfill
		\begin{minipage}{0.45\textwidth}
			\centering
			\includegraphics[width=\textwidth]{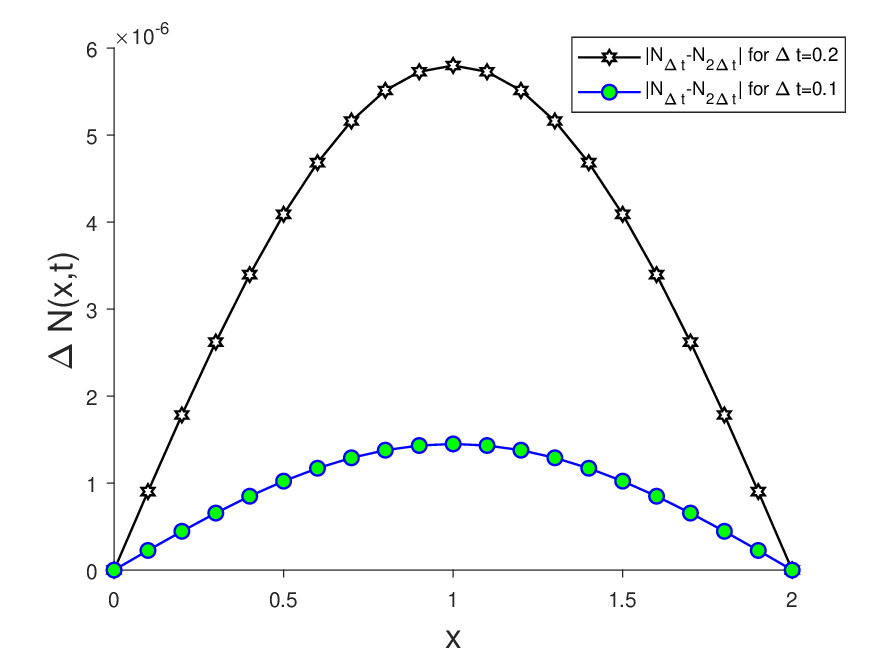}
			\subcaption{}
		\end{minipage}		
		\caption{Absolute error of $M(x,t)$ and $N(x,t)$ from equation (\ref{eq:4.1}) for different time increment at time $t=10$.}
\label{fig:1nerr}
	\end{figure}

\noindent\emph{\textbf{Test Problem 1: }}Let us consider the system of the parabolic equations from the study of Manaa \textit{et. al.}\cite{r37}
\begin{align}
\setstretch{2.0}
\left.
\begin{array}{lcl}
\displaystyle \frac{\partial M}{\partial t}=\varepsilon_1 \frac{\partial^2 M}{\partial x^2} + f(M,N) -(p+q)M\\
\displaystyle \frac{\partial N}{\partial t}=\varepsilon_2\frac{\partial^2 N}{\partial x^2}- f(M,N) +p(1-N)\\
\end{array}
\right\}
\label{eq:4.1}
\end{align}
where $f(M,N)=M^2N$ and $\displaystyle x\in [a,b], t\ge 0$. The boundary conditions and the initial conditions are considered as:
\begin{align}
\setstretch{1.0}
\left.
    \begin{array}{lcl}
    \displaystyle M(a,t)=M(b,t)=0\\
    \displaystyle N(a,t)=N(b,t)=1\\
    \end{array}
    \right\}
    \label{eq:4.2}
\end{align}
\noindent and
  \begin{align}
  \setstretch{1.0}
    \left.
    \begin{array}{lcl}
    \displaystyle M(x,0)=0.01 sin(\pi (x-b)/(b-a))\\
    \displaystyle N(x,0)=1-0.12 sin(\pi (x-b)/(b-a))\\
    \end{array}
    \right\}
    \label{eq:4.3}
\end{align}
\noindent The domain of the model is $[a, b]$. The values of the parameters are taken as $ a=0, b=2, \varepsilon_1=\varepsilon_2=0.01$, $p=0.09$, and $q=-0.004$.\\
\noindent Here, to obtain the numerical approximation, the effect of boundary conditions is insignificant because all terms of $B_j(x)$ are zero at the boundary points. We have employed the modified Galerkin method to the system of nonlinear partial differential equations (\ref{eq:4.1}) and therefore obtained the system of ordinary differential equations with respect to $t$. In this stage, we have used the $\alpha$ family of approximation in order to convert the system into recurrent relations and then we applied \textit{Picard iterative procedure}. To find the initial guess of the given system, we have applied the weighted residual procedure on the initial conditions (\ref{eq:4.2}).\\
\noindent Tables (\ref{tab:4.1}) and (\ref{tab:4.11}) provide the numerical results of concentrations $M (x,t)$ and $N(x,t)$ for various values of $x$. For computation, we have taken $\Delta t=0.1$. The numerical approximations are derived at time levels $t=1$ and $t=2$. 
\begin{table}[ht]
\centering
\caption{Numerical results of concentrations $M (x,t)$ at different time levels with $\Delta t=0.1$ and first 7 modified Bernstein polynomials.}
\setstretch{1.0}
   \resizebox{0.9\textwidth}{!}{
   \begin{tabular}{ c | c | c | c | c }
   \hline
   \multirow{2}{*}{\textbf{\textit{x}}}   &\multicolumn{2}{c}{$t=1$}&\multicolumn{2}{|c}{$t=2$}\\
  \cline{2-5}
  &\textbf{Present Method}&\textbf{Reference \cite{r37}}&\textbf{Present Method}&\textbf{Reference \cite{r37}}\\
 \hline
 0.0&0.00&0.00&0.00&0.00\\
 \hline
 0.1&-0.00139&-0.0163&-0.00124&-0.0139\\
 \hline
 0.2&-0.00275&-0.0317&-0.00243&-0.0271\\
 \hline
 0.3&-0.00402&-0.0461&-0.00356&-0.0390\\
 \hline
 0.4&-0.00520&-0.0589&-0.00459&-0.0495\\
 \hline
 0.5&-0.00624&-0.0701&-0.00550&-0.0585\\
 \hline
 0.6&-0.00712&-0.0793&-0.00627&-0.0659\\
 \hline
 0.7&-0.00782&-0.0867&-0.00688&-0.0716\\
 \hline
 0.8&-0.00834&-0.0920&-0.00733&-0.0757\\
 \hline
 0.9&-0.00866&-0.0952&-0.00760&-0.0782\\
 \hline
 1.0&-0.00877&-0.0962&-0.00769&-0.0790\\
 \hline
 1.1&-0.00866&-0.0952&-0.00760&-0.0782\\
 \hline
 1.2&-0.00835&-0.0920&-0.00733&-0.0757\\
 \hline
 1.3&-0.00782&-0.0867&-0.00688&-0.0716\\
 \hline
 1.4&-0.00712&-0.0793&-0.00626&-0.0659\\
 \hline
 1.5&-0.00623&-0.0701&-0.00550&-0.0585\\
 \hline
 1.6&-0.00520&-0.0589&-0.00459&-0.0495\\
 \hline
1.7&-0.00402&-0.0461&-0.00356&-0.0390\\
\hline
1.8&-0.00275&-0.0318&-0.00243&-0.0271\\
\hline
1.9&-0.00139&-0.0162&-0.00124&-0.0139\\
\hline
2.0&0.00&0.0003&0.00&-0.0000\\
\hline
\end{tabular}
}
	\label{tab:4.1}
\end{table}
  \FloatBarrier
\noindent Throughout these tables, we have compared the results which we have obtained with the numerical approximations that have already been published in other well-known literature. The table demonstrates that our outcomes are reasonably comparable to those that have been published. It validates the accuracy of our approach to approximating the reaction-diffusion system numerically.
\begin{table}[ht]
\centering
  \caption{Numerical results of concentrations $N (x,t)$ at different time levels with $\Delta t=0.1$ and first 7 modified Bernstein polynomials.}
\setstretch{1.0}
   \resizebox{0.9\textwidth}{!}{
   \begin{tabular}{c | c | c | c | c}
   \hline
   \multirow{2}{*}{\textbf{\textit{x}}}   &\multicolumn{2}{c}{$t=1$}&\multicolumn{2}{|c}{$t=2$}\\
  \cline{2-5}
  &\textbf{Present Method}&\textbf{Reference \cite{r37}}&\textbf{Present Method}&\textbf{Reference \cite{r37}}\\
 \hline
 0.0&1.0000&1.0000&1.0000&1.0000\\
 \hline
 0.1&1.01673&1.0014&1.01491&1.0014\\
 \hline
 0.2&1.03305&1.0027&1.02946&1.0026\\
 \hline
 0.3&1.04855&1.0039&1.04327&1.0035\\
 \hline
 0.4&1.06285&1.0048&1.05601&1.0043\\
 \hline
 0.5&1.0756&1.0055&1.06736&1.0048\\
 \hline
 0.6&1.08650&1.0061&1.07706&1.0052\\
 \hline
 0.7&1.09525&1.0065&1.08485&1.0054\\
 \hline
 0.8&1.10167&1.0068&1.09056&1.0056\\
 \hline
 0.9&1.10558&1.0070&1.09404&1.0056\\
 \hline
 1.0&1.10689&1.0070&1.09521&1.0057\\
 \hline
 1.1&1.10558&1.0070&1.09404&1.0056\\
 \hline
 1.2&1.10167&1.0068&1.09056&1.0056\\
 \hline
 1.3&1.09525&1.0065&1.08485&1.0054\\
 \hline
 1.4&-1.08650&1.0061&1.07705&1.0052\\
 \hline
 1.5&1.07561&1.0055&1.06736&1.0048\\
 \hline
 1.6&1.06285&1.0048&1.05061&1.0043\\
 \hline
1.7&1.04855&1.0039&1.04327&1.0035\\
\hline
1.8&1.03305&1.0027&1.02945&1.0026\\
\hline
1.9&1.01673&1.0014&1.01491&1.0014\\
\hline
2.0&1.0000&1.0000&1.0000&1.0000\\
\hline
\end{tabular}
}
	\label{tab:4.11}
\end{table}
  \FloatBarrier
\noindent The approximate results $M(x, t)$ and $N(x, t)$ of Equation  (\ref{eq:4.1}) are presented in the following figure (\ref{fig:1m}). 
\begin{figure}[ht!]
		\centering
		\begin{minipage}{0.45\textwidth}
			\centering
			\includegraphics[width=\textwidth]{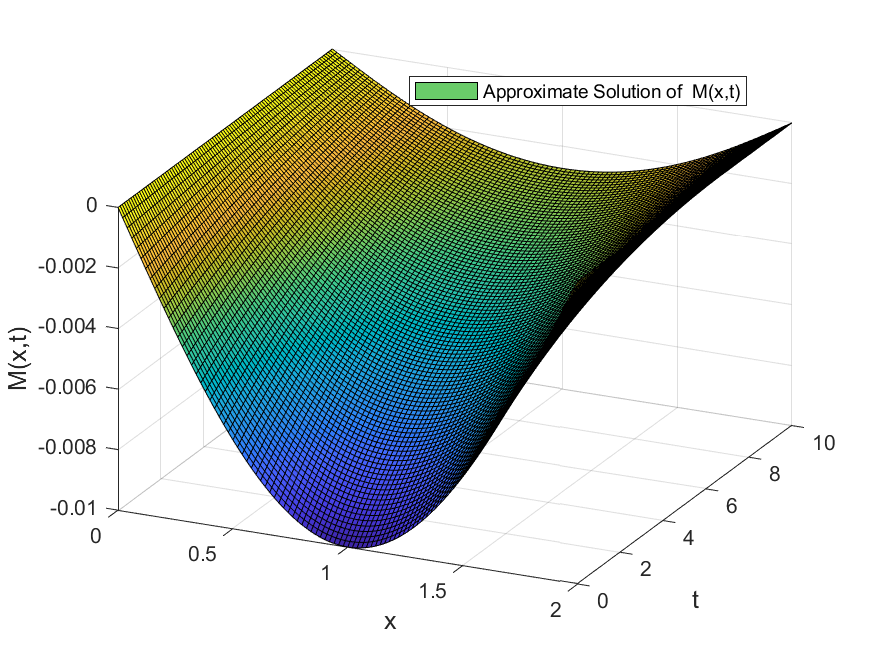}
			\subcaption{}
		\end{minipage}\hfill
		\begin{minipage}{0.45\textwidth}
			\centering
			\includegraphics[width=\textwidth]{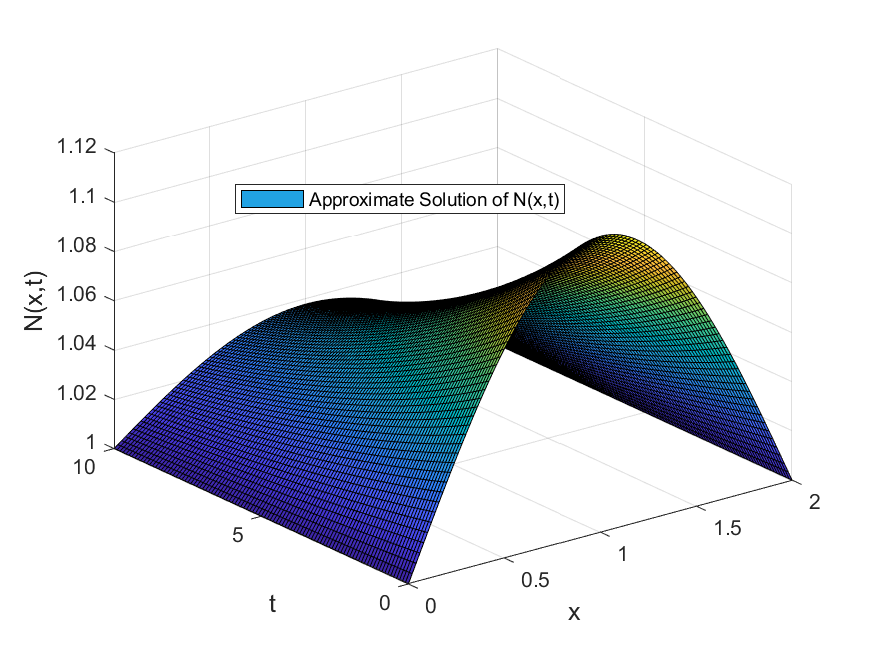}
			\subcaption{}
		\end{minipage}		
		\caption{Approximate solution of $M(x,t)$ and $N(x,t)$ of (\ref{eq:4.1}) by using the present method for $(x,t)\in [0,2]\times [0,2]$}
		\label{fig:1m}
	\end{figure}
 \FloatBarrier
\noindent In Figure (\ref{fig:1m}) we have employed a three-dimensional graphical depiction of approximate solutions of $M(x,t)$ and $N(x,t)$ at different time levels for better understanding. The graphical representations agree with the results that we have obtained in the tables. Eventually, it makes sense clearly that the method is more applicable to solving such nonlinear parabolic PDE systems.\\
In Figure (\ref{fig:1nerr}), we have presented the error graph of $M(x,t)$ and $N(x,t)$ at time $t=10$, where the absolute errors are computed between two different time increments, say $\Delta t=0.2$, $\Delta t=0.4$ and $\Delta t=0.1$, $\Delta t=0.2$.
 \FloatBarrier

\noindent The $L_2$ norm and  $L_\infty$ norm, are presented in Table (\ref{tab:prob1Mt2}), which shows that the comparative errors are reduced significantly according to the reduction of the size of the time increments.
\begin{table}[ht]
\caption{The $L_2$ and $ L_\infty$ norm  at $t=10$ for $M(x,t)$ of equation (\ref{eq:4.1}).}
\medskip
\centering\renewcommand{\arraystretch}{1.0}
\begin{tabular}{c|c|c|c|c}
\hline
   \multirow{2}{*}{$\Delta t$}   &\multicolumn{2}{c}{$M(x,t)$}&\multicolumn{2}{|c}{$N(x,t)$}\\
  \cline{2-5}
  &\textbf{$L_2$ norm}&\textbf{$L_\infty$ norm}&\textbf{$L_2$ norm}&\textbf{$L_\infty$ norm}\\
\hline

0.40 & -  & -&-&-\\ \hline
0.20 & 0.00000138&0.00000060&0.00001359&0.00000580\\ \hline
0.10 &  0.00000035& 0.00000015&0.00000340&0.00000145\\\hline
\end{tabular}

\label{tab:prob1Mt2}
\end{table}
\FloatBarrier

\noindent \emph{\textbf{Test Problem 2: }} The Gray-Scott Model is one of the most important models whose wave formations are similar to many waves formed in real life such as butterfly wings, gesticulation, damping, turning patterns, embryos, multiple spots, and so on \cite{r12,r38}. Let us consider the following model,
\begin{align}
\setstretch{2.0}
\left.
\begin{array}{lcl}
\displaystyle \frac{\partial M}{\partial t}=\varepsilon_1 \frac{\partial^2 M}{\partial x^2} - f(M,N) +p(1-M)\\
\displaystyle \frac{\partial N}{\partial t}=\varepsilon_2\frac{\partial^2 N}{\partial x^2}+ f(M,N) -(p+q)N\\
\end{array}
\right\}
\label{eq:4.5}
\end{align}
where $f(M,N)=MN^2$. The boundary conditions and the initial conditions are considered as follows:
\begin{align}
\setstretch{1.0}
    \left.
    \begin{array}{lcl}
    \displaystyle M(-50,t)=M(50,t)=1\\
    \displaystyle N(-50,t)=N(50,t)=0\\
    \end{array}
    \right\}
    \label{eq:4.6}
\end{align}
\noindent and
  \begin{align}
  \setstretch{1.0}
    \left.
    \begin{array}{lcl}
    \displaystyle M(x,0)=1-0.5 sin^{100}(\pi (x-50)/100)\\
    \displaystyle N(x,0)=0.25 sin^{100}(\pi (x-50)/100)\\
    \end{array}
    \right\}
    \label{eq:4.7}
\end{align}
\noindent The domain of the model is $[-50, 50]$. The values of the parameters are taken as 
\begin{center}
$\varepsilon_1=1$, $\varepsilon_2=0.01$, $p=0.01$, $q=0.12$
\end{center}
\noindent Here for computational purposes, we have used $7$ modified Bernstein polynomials. By applying the modified Galerkin method, we have used the backward difference method to transform the system of ordinary differential equations into the recurrent relations which is therefore solved by \textit{Picard} iterative procedure. \noindent Numerical data of $M(x, t)$ and $N(x, t)$ of (\ref{eq:4.1}) are also presented in tabulated form in the following table at different time steps.
\begin{table}[ht]
\caption{Numerical results of concentrations $M(x,t)$ and $N(x,t)$ ate different time levels with $\Delta t=0.1$ and first 7 modified Bernstein polynomials.}
\centering
   \resizebox{0.8\textwidth}{!}{
   \begin{tabular}{ c | c | c | c | c | c | c }
   \hline
   \multirow{2}{*}{\textbf{\textit{x}}} &\multicolumn{3}{c|}{\textbf{$M(x,t)$}}&\multicolumn{3}{|c}{\textbf{$N(x,t)$}}\\
  \cline{2-7}
&\textbf{\textit{$t=1$}}&\textbf{\textit{$t=10$}}&\textbf{\textit{$t=20$}}&\textbf{\textit{$t=1$}}&\textbf{\textit{$t=10$}}&\textbf{\textit{$t=20$}}\\
  \hline
 -50.0&1.0000&1.0000&1.0000&0.0000&0.0000&0.0000\\
  \hline
  -40.0&0.9725&0.9746& 0.9767&0.0140&0.0164&0.0209\\
  \hline
  -30.0&1.0155&1.0201&1.0192&-0.0071&-0.0072&-0.0094\\
  \hline
  -20.0&1.0288&1.0056&0.9866&-0.0144&-0.0073&-0.0061\\
  \hline
  -10.0&0.8770&0.8515& 0.8284&0.0623&0.0850&0.1148\\
  \hline
  0.0&0.7753&0.7563&0.7355&0.1139&0.1447&0.1922\\
  \hline
  10.0&0.8770&0.8515&0.8284&0.0623&0.0850&0.1148\\
  \hline
  20.0&1.0288&1.0056&0.9866&-0.0144&-0.0073&-0.0061\\
  \hline
  30.0&1.0155&1.0201&1.0192&-0.0071&-0.0072&-0.0094\\
  \hline
  40.0&0.9725&0.9746&0.9767&0.0140&0.0164&0.0209\\
  \hline
  50.0&1.0000&1.0000&1.0000&0.0000&0.0000&0.0000\\
 \hline
 \end{tabular}
  }
  
  \label{tab:4.2}
  \end{table}
  \FloatBarrier
\noindent The table shows that the numerical values of concentrations $M$ and $N$ change very slowly with varying values of $x$. It happens in every time step.

\begin{figure}[ht!]
		\centering
		\begin{minipage}{0.45\textwidth}
			\centering
			\includegraphics[width=\textwidth]{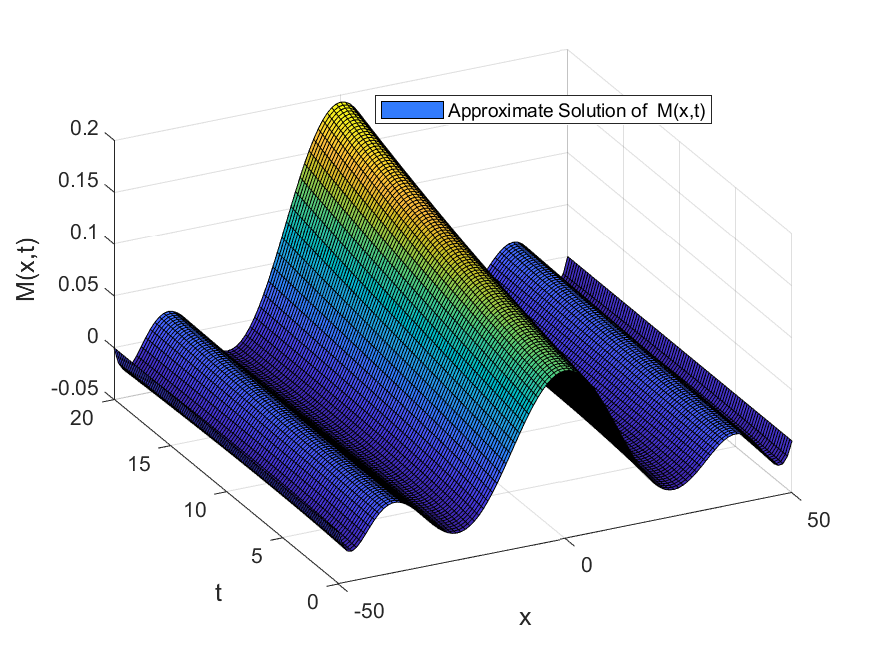}
			\subcaption{}
		\end{minipage}\hfill
		\begin{minipage}{0.45\textwidth}
			\centering
			\includegraphics[width=\textwidth]{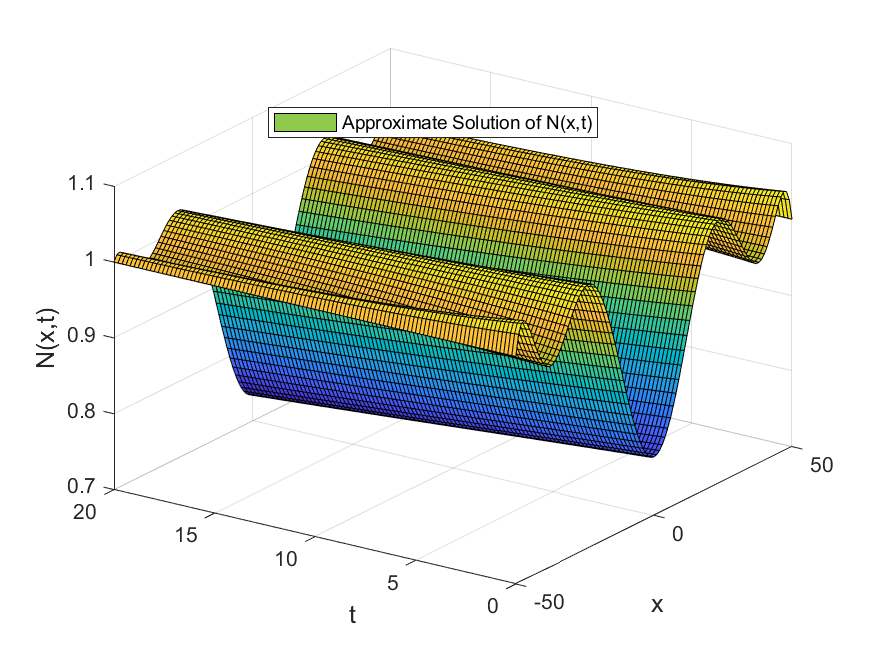}
			\subcaption{}
		\end{minipage}		
		\caption{Absolute errors of $M(x,t)$ and $N(x,t)$ of (\ref{eq:4.5}) by using the present method for $(x,t)\in [-50,50]\times [0,20]$}
		\label{fig:2m}
	\end{figure}

\noindent The results obtained by applying our proposed scheme are presented in Figure (\ref{fig:2m}).
 \FloatBarrier
\noindent Figure (\ref{fig:2m}) is deployed to provide pictorial representations of the numerical concentrations  $M$ and $N$ at different time levels. The results that are obtained in the table are shown graphically. The graphs are obtained for different time levels. The graphical presentation shows that the changes in concentrations are sufficiently small for different time levels.\\

\begin{figure}[ht!]
		\centering
		\begin{minipage}{0.45\textwidth}
			\centering
			\includegraphics[width=\textwidth]{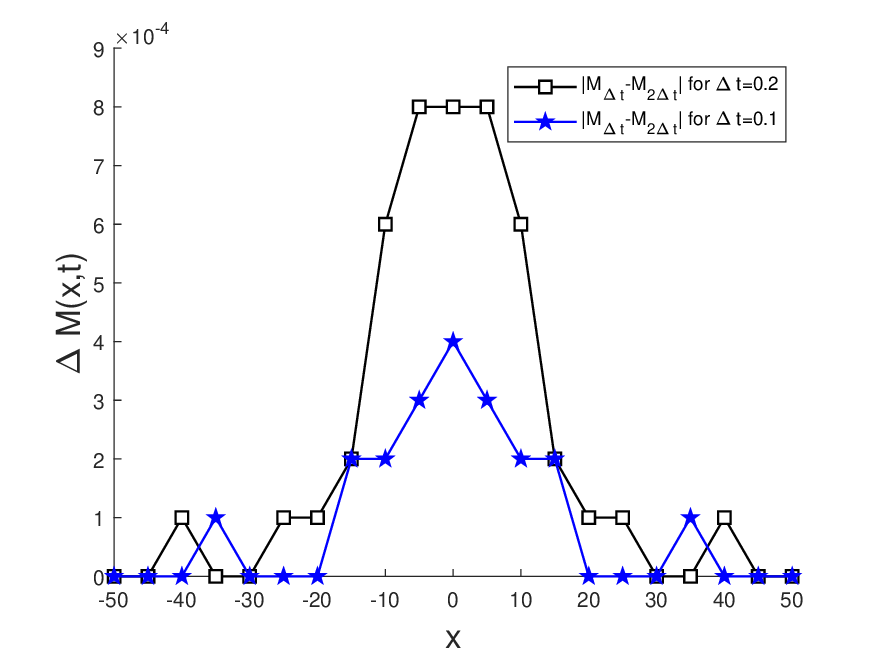}
			\subcaption{}
		\end{minipage}\hfill
		\begin{minipage}{0.45\textwidth}
			\centering
			\includegraphics[width=\textwidth]{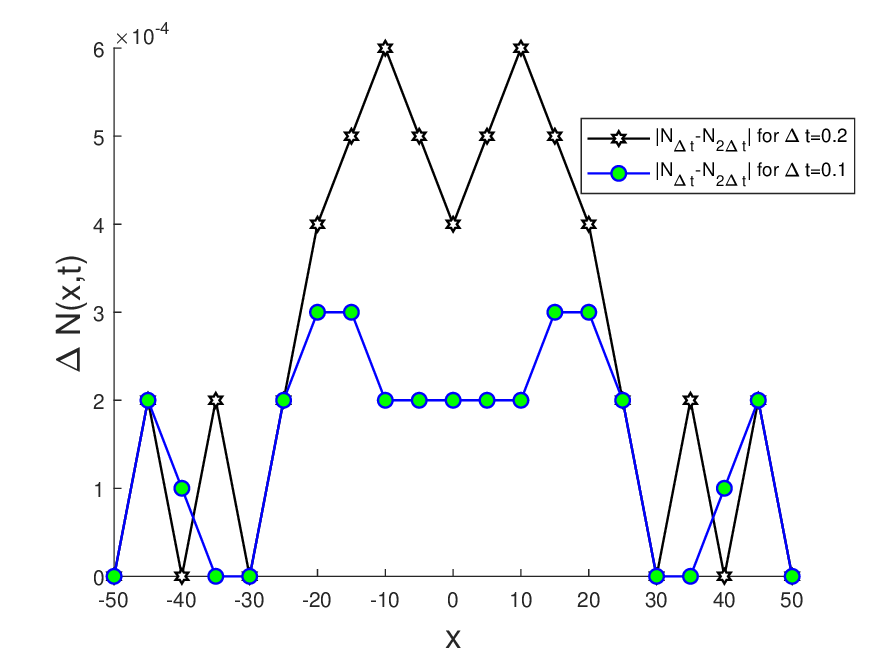}
			\subcaption{}
		\end{minipage}		
		\caption{Approximate solution of $M(x,t)$ and $N(x,t)$ of (\ref{eq:4.1}) by using the present method for $(x,t)\in [0,2]\times [0,2]$}
		\label{fig:2merr}
	\end{figure}

\noindent The $L_2$ norm, and  $L_\infty$ norms, are presented in table (\ref{tab:prob2Mt2}), which shows that the comparative errors are reduced significantly according to the reduction of the size of the time increments. However, the order of convergences increased noticeably along with the reduction of the time length.
\begin{table}[!ht]
\caption{The $L_2$ and $ L_\infty$ norm  at $t=10$ for $M(x,t)$ of equation (\ref{eq:4.5}).}
\centering\renewcommand{\arraystretch}{1.0}
\begin{tabular}{c|c|c|c|c}
\hline
   \multirow{2}{*}{$\Delta t$}   &\multicolumn{2}{c}{$M(x,t)$}&\multicolumn{2}{|c}{$N(x,t)$}\\
  \cline{2-5}
  &\textbf{$L_2$ norm}&\textbf{$L_\infty$ norm}&\textbf{$L_2$ norm}&\textbf{$L_\infty$ norm}\\
\hline
0.40  & - & -&-&- \\ \hline
0.20  & 0.00110596&0.00080091&0.00105161&0.00051358\\ \hline
0.10  & 0.00055524& 0.00040251&0.00052415& 0.00025600\\\hline
\end{tabular}
\label{tab:prob2Mt2}
\end{table}
\FloatBarrier
\noindent In Figure (\ref{fig:2merr}), we have presented the error graph of $M(x,t)$ and $N(x,t)$ at time $t=10$, where the absolute errors are computed between two different time increments say $\Delta t=0.2$, $\Delta t=0.4$ and $\Delta t=0.1$, $\Delta t=0.2$. 
 \FloatBarrier
\section*{Conclusion}
\noindent This research study has provided numerical approximations of nonlinear reaction-diffusion systems with specified boundary and initial conditions through the employment of the modified Galerkin method. To generate the trial solution, modified Bernstein Polynomials have been used. The simplification of the weighted residual leads to a system of ordinary differential equations which is then transformed into the recurrent relation by applying the backward difference formula. At this stage, we have used Picard's iterative procedure to approximate the trial solution. After successful derivation, we applied our proposed method to several models in order to test their applicability and effectiveness. We have solved and displayed the results both numerically and graphically. From those figures and numerical results, it is indisputable that our proposed method is an unconditionally stable, efficient, highly modular, and easily expandable method that can be applied to any type of system of nonlinear parabolic partial differential equations regardless of the type of the boundary conditions, type of non-linearity of the functions, coefficients are constants or function of independent variables.
\section*{Acknowledgement}
The authors acknowledge that the research was supported and funded by Dhaka University research grant under UGC, Bangladesh.

        \end{document}